\documentclass[review]{elsarticle}
\usepackage{lineno}
\usepackage{hyperref}
\usepackage{bbold}
\usepackage{verbatim}
\usepackage{url}
\usepackage{graphicx}
\usepackage{float}
\usepackage{mathrsfs}
\usepackage{epstopdf}
\usepackage[all,color]{xy}
\usepackage{cancel}
\usepackage{MnSymbol}
\usepackage{tikz}
\usepackage[margin=1in]{geometry}

\usetikzlibrary{arrows}
\modulolinenumbers[5]
\allowdisplaybreaks

\newtheorem{theorem}{Theorem}

\newtheorem{corollary}[theorem]{Corollary}
\newtheorem{example}[theorem]{Example}
\newtheorem{lemma}[theorem]{Lemma}
\newtheorem{proposition}[theorem]{Proposition}

\newproof{proof}{Proof}
\numberwithin{equation}{subsection}
\numberwithin{theorem}{subsection}

\newcommand{\ignore}[1]{}

\DeclareMathOperator{\Hom}{Hom}

\begin{document}

\begin{frontmatter}

\title{Oriented Hypergraphs: Balanceability}

\author[add2]{Lucas J. Rusnak}\corref{mycorrespondingauthor}\ead{Lucas.Rusnak@txstate.edu}
\author[add1]{Selena Li}
\author[add1]{Brian Xu}
\author[add1]{Eric Yan}
\author[add1]{Shirley Zhu}

\address[add2]{Department of Mathematics, Texas State University, San Marcos, TX 78666, USA}

\address[add1]{Mathworks, Texas State University, San Marcos, TX 78666, USA}

%\fntext[fn2]{Portions of these results appear in 2019 Honors thesis.}

\cortext[mycorrespondingauthor]{Corresponding author}

\begin{abstract}
An oriented hypergraph is an oriented incidence structure that extends the concepts of signed graphs, balanced hypergraphs, and balanced matrices. We introduce hypergraphic structures and techniques that generalize the circuit classification of the signed graphic frame matroid to any oriented hypergraphic incidence matrix via its locally-signed-graphic substructure. To achieve this, Camion's algorithm is applied to oriented hypergraphs to provide a generalization of reorientation sets and frustration that is only well-defined on balanceable oriented hypergraphs. A simple partial characterization of unbalanceable circuits extends the applications to representable matroids demonstrating that the difference between the Fano and non-Fano matroids is one of balance.
\end{abstract}

\begin{keyword}
Oriented hypergraph \sep balanced hypergraph \sep balanced matrix \sep balancing sets \sep signed graph.
\MSC[2010] 05C75 \sep 05C65 \sep 05C22 \sep 05C50 \sep 05B35 
\end{keyword}

\end{frontmatter}

%%%%%%%%
% use these commands for typesetting doi and arXiv references in the bibliography
% if needed include a line break (\\) at an appropriate place in the title
% \date{\dateline{submission date}{acceptance date}\\
% \small Mathematics Subject Classifications: comma separated list of
% MSC codes available from http://www.ams.org/mathscinet/freeTools.html}

%% main text

\section{Introduction}

An oriented hypergraph is a signed incidence structure where each unique vertex-edge incidence is given a label of $+1$ or $-1$ and each adjacency is signed as the negative of the product of its incidences. Matrices with commensurable entries can be represented as an oriented hypergraph via multiple incidences of unit weight, and may subsequently be studied via their locally-signed-graphic substructure. Oriented hypergraphs provide a way to merge and generalize the study of balanced hypergraphs \cite{Berge1, CO1},  balanced $\{0,  \pm1 \}$ matrices \cite{DBM, TrAlpha,TrLog}, and algebraic graph theory \cite{AH1, OH1}. Spectral properties of oriented hypergraphs have been studied in \cite{Reff6, Reff2}, while various characteristic polynomials of the adjacency and Laplacian matrices of oriented hypergraphs were classified in \cite{OHSachs, OHMTT}, providing a unifying generalization of matrix-tree-type Theorems and Sachs-type Theorems, as well as a solution to the maximum permanent that can be refined to solve the maximum determinant. 

It was shown in \cite{IH1} that incidence hypergraphs provide the central point to study combinatorial matrix theory, where the incidence matrix and the bipartite representation graph are natural Kan extensions of logical functor in the category of incidence hypergraphs --- this approach solves the characterization of graph exponentials and shows they are $\Hom$'s in the category of incidence hypergraphs. The subobject classifier of the topos was used to provide a characterization of the all-minors characteristic polynomial for integer matrices via subhypergraphic families in the injective envelope in \cite{IH2}.

The focus of this paper is to continue the hypergraphic structural characterization of the circuits of vector matroids with commensurable entries introduced in \cite{OH1}, where the structure of oriented hypergraphs was broken into three main categories -- balanced, balanceable, and unbalanceable -- and the balanced circuits were characterized. We provide a strengthening of these results by extending the circuit characterization to include balanceable circuits; hence, contain the circuit characterization of the signed graphic frame matroid by Zaslavsky in \cite{SG}. The characterization is accomplished by adapting techniques from balanced matrices and re-interpreting Camion's algorithm \cite{Camion,BM} as a method for incidence re-orientation that shifts the notion of a balancing set and the frustration index to the underlying incidence structure. These balancing sets of incidences are used to introduce an arterial connection (hypergraphic path) called \emph{shuntings} to produce circuits. 

Finding a characterization of the unbalanceable circuits would complete the circuit classification for these vector matroids using oriented hypergraphic families, and hopefully lend to new graph-like techniques that may be used to study representable matroids. The techniques of Camion do not extend to unbalanceable hypergraphs and the frustration index (as determined by switching) is no longer well-defined. However, a simple characterization of the unbalanceable circuits arising from a single minimal cross-theta is provided along with how their negative circle structure forces them to vanish modulo $k$. Finally, the difference between the Fano and non-Fano matroids is shown to be a balance property. Combined with \cite{IH1} it opens the door for a modified version of graph theoretic techniques to be applied to representable matroids.

\section{Background}
\subsection{Oriented Hypergraph Basics}

The definitions in this section are condensed from \cite{OH1} and updated following the work done in \cite{OHSachs, IH1}. An \emph{oriented hypergraph} is a quintuple $(V,E,I,\iota,\sigma)$ consisting of a set of vertices $V$, a set of edges $E$, a set of incidences $I$, an incidence function $\iota:I\rightarrow V\times E$, and an orientation function $\sigma:I\rightarrow \{+1,-1\}$. A value of $+1$ is indicated by an arrow at incidence $i$ entering the vertex, while a value of $-1$ is indicated by an arrow at incidence $i$ exiting the vertex. The \emph{incidence dual} of an oriented hypergraph $G$ is the oriented hypergraph $G^*$ where the vertex set and edge set are reversed. An oriented hypergraph in which each edge is assigned exactly two incidences is called a \emph{bidirected graph}. A bidirected graph in which every edge/adjacency is positive is regarded as an orientation of an ordinary graph as they have indistinguishable incidence matrices; see \cite{MR0267898, SG, OSG} for bidirected graphs as orientations of signed graphs. The \emph{incidence matrix} of an oriented hypergraph $G$ is the $V \times E$ matrix $\mathbf{H}_{G}$ where the $(v,e)$-entry is the sum of $\sigma(i)$ for each $i \in I$ such that $\iota(i)=(v,e)$. The \emph{bipartite representation graph of $G$} is the bipartite graph $\Gamma$ where $V(\Gamma)= V \cup E$ and $E(\Gamma)= I$.

\begin{figure}[H]
    \centering
    \includegraphics{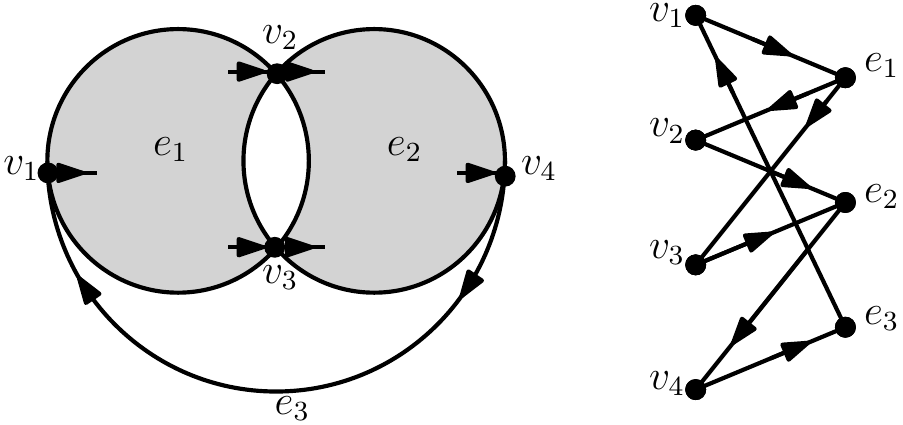}
    \caption{An oriented hypergraph and its bipartite representation.}
    \label{fig:OH}
\end{figure}

Many of the oriented hypergraphic definitions coincide with either a locally-signed-graphic embedding \cite{OHHar, AH1, OH1} or by translating the graphic definitions from the bipartite incidence graph $\Gamma$ back to $G$ via the corresponding atomic geometric morphism \cite{IH1}.

A \emph{directed path of length $n/2$} is a non-repeating sequence 
\begin{equation*}
\overrightarrow{P}_{n/2}=(a_{0},i_{1},a_{1},i_{2},a_{2},i_{3},a_{3},...,a_{n-1},i_{n},a_{n})
\end{equation*}
of vertices, edges, and incidences, where $\{a_k\}$ is an alternating sequence of vertices and edges, and $i_{k}$ is an incidence between $a_{k-1}$ and $a_{k}$. A \emph{circle} is a closed directed path. A \emph{directed adjacency of $G$} is an incidence-monic map of $\overrightarrow{P}_{1}$ into $G$. The \emph{sign of a path} $P$ is 
\begin{equation*}
sgn(P)=(-1)^{\lfloor n/2\rfloor }\prod_{k=1}^{n}\sigma (i_{k})\text{,}
\end{equation*}
which is equivalent to taking the product of the signed adjacencies if $P$ is a vertex-path. 

An oriented hypergraph is \emph{inseparable} if every pair of incidences is contained in a circle. For convention, a $1$-edge is not inseparable, but a $0$-edge is inseparable. A \emph{flower} is a minimally inseparable oriented hypergraph; the oriented hypergraph in Figure \ref{fig:OH} is an example of a flower. A flower is the hypergraphic generalization of a circle in a bidirected graph.

\begin{proposition}[\cite{OH1}, Prop. 4.1.2]
\label{onlyflowers}
$F$ is a flower of a signed graph if, and only if, $F$ is a circle or a loose edge.
\end{proposition}

A monovalent vertex $v$ is a \emph{thorn} of an oriented hypergraph $G$ if some circle of $G$ contains the edge incident to $v$. A \emph{pseudo-flower} is an oriented hypergraph containing one or more thorns, where weak-deletion of all thorns results in a flower, called the \emph{flower-part} --- that is, set deletion of the thorn-vertices and removal of the incidence, but leaving the edge untouched. The \emph{cyclomatic number of $G$} is 
\begin{align*}
\varphi=|I|-(|V|+|E|)+c,
\end{align*}
where $c$ is the number of connected components of $G$; this is equivalent to the cyclomatic number of $\Gamma$ as there is one-to-one correspondence between their circles.

A \emph{subdivision} of a $k$-edge $e$ with incidence set $I(e)=\{i_1, i_2, \ldots i_k\}$ replaces $e$ with two new edges $e_1$ and $e_2$, and introduces two new incidences $j_1$ and $j_2$ and a new vertex $w$ incident to $j_1$ and $j_2$ such that $j_1 \in I(e_1)$, $j_2 \in I(e_2)$, $I(e_1) \cap I(e_2) = \emptyset$, and $(I(e_1) \setminus j_1)  \cup (I(e_2) \setminus j_2) = I(e)$. An edge subdivision is \emph{incompatible} if $\sigma(j_1)\sigma(j_2)=+1$, and \emph{compatible} if $\sigma(j_1)\sigma(j_2)=-1$ --- this is equivalent to the sign of the co-adjacency being negative (incompatible) or positive (compatible).
An \emph{artery} is either a single vertex, or a subdivision of a $k$-edge ($k \geq 2$). The degree-$2$ vertices of an artery are called \emph{internal} vertices, while the non-degree-$2$ vertices are \emph{external}. An artery is the hypergraphic generalization of a path. 

\begin{proposition}[\cite{OH1}, Prop. 4.2.3]
\label{onlyartery}
$A$ is an artery of a signed graph if, and only if, $A$ is a path.
\end{proposition}

An \emph{arterial connection of pseudo-flowers by thorns} is a collection of pseudo-flowers connected to each other via arteries so that no new circles are created --- two pseudo-flowers are allowed to be connected by a single vertex-artery that is a thorn common to both.

\subsection{Balance, Thetas, Circuits, and Frustration}

An oriented hypergraph is \emph{balanced} if the sign of each circle is positive; an oriented hypergraph is \emph{balanceable} if there exists a balanced orientation; and an oriented hypergraph is \emph{unbalanceable} if it is not balanceable. A \emph{cross-theta} is a subgraph of an oriented hypergraph that consists of three internally disjoint paths of half-integer length; equivalently, three internally disjoint paths of odd length in the bipartite representation graph. Analogous definitions are used for a \emph{vertex-theta} and an \emph{edge-theta}.

\begin{figure}[H]
    \centering
    \includegraphics{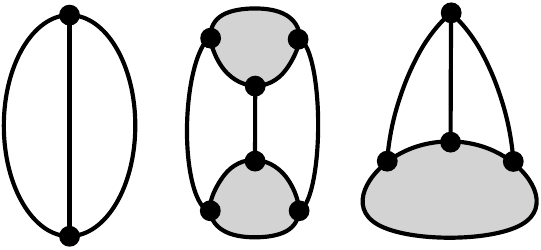}
    \caption{A vertex-, edge-, and cross-theta, respectively.}
    \label{fig:thetas}
\end{figure}

The following theorems emphasize the importance of cross-thetas. 

\begin{theorem}[\cite{OH1}, Theorem 5.3.8.]
\label{ohcross}
If a flower contains a vertex of degree $\geq$ 3, then it contains a cross-theta.
\end{theorem}

\begin{theorem}[\cite{OH1}, Prop. 6.2.2]
\label{crossthetabad}
An oriented hypergraph $G$ is balanceable if, and only if, it does not contain a cross-theta.
\end{theorem}

Additionally, balanceable flowers have the following strong condition on their bipartite representation graph.

\begin{theorem}[\cite{OH1}, Theorem 5.4.2]
\label{ears}
If $F$ is a cross-theta-free flower, then every ear decomposition of $\Gamma_F$ can be regarded as consisting of only edge-paths --- that is, paths that start and end on the edge-vertex side of $\Gamma_F$.
\end{theorem}

The characterization of the circuits of the graphic and signed graphic frame matroid (see \cite{AGT,SG}) have a simple reinterpretation in terms of oriented hypergraphs.

\begin{theorem}
The circuits of the graphic matroid are flowers.
\end{theorem}
\begin{proof}
The only circuits of the graphic matroid are graphic circles. The result follows from Proposition \ref{onlyflowers}. 
\qed \end{proof}

\begin{theorem}
The circuits of the signed graphic frame matroid are balanced flowers, or arterially-connected unbalanced flowers.
\end{theorem}
\begin{proof}
The circuits of the signed graphic frame matroid were characterized by Zaslavsky \cite{SG} and are either positive circles, or two negative circles connected by a path (with $1$-edges regarded as negative loops). Since $1$-edges are pseudo-flowers, the result follows from Propositions \ref{onlyflowers} and \ref{onlyartery}.
\qed \end{proof}

The balanced circuits of any oriented hypergraph have also been characterized, but require a bit more background. The \emph{balanced subdivision} of an edge is any subdivision in which the signs of corresponding circles do not change (i.e. the only incompatible subdivisions involve adjacencies not in any circle). The incidence inverse operation to subdivision is \emph{$2$-vertex-contraction}, which is equivalent to the signed graphic contraction of the corresponding $2$-edge in the incidence dual --- this is an inverse on the incidence structure and on compatible subdivision.

\begin{lemma}[\cite{OH1}, Lemmas 3.1.5 \& 3.2.3]
\label{balsubgood}
Let $H$ be an edge-induced subhypergraph of $G$, and let $H'$ be obtained by balanced subdivision. $H$ is a circuit if, and only if $H'$ is a circuit.
\end{lemma}

\begin{theorem}[\cite{OH1}, Prop. 6.2.7]
\label{BalMD}
The balanced circuits of an oriented hypergraph are balanced flowers, or have a balanced subdivision that is an 
arterial connection of pseudo-flowers by thorns.
\end{theorem}

We obtain a characterization of the balanceable circuits of incidence matrices associated to an oriented hypergraph. Since signed graphs are balanceable oriented hypergraphs, the remainder of Zaslavsky's circuit characterization in \cite{SG} is a corollary of this result. This is obtained by translating the concept of frustration to oriented hypergraphs.

\subsection{Frustration in Oriented Hypergraphs}

Harary introduced in \cite{Har1} the frustration index of a signed graph as the smallest number of edges whose deletion (equivalently, negation) results in a balanced signed graph. Such a set of edges is called a \emph{balancing set}.

A \emph{switching function} on a signed graph is any function $\varsigma : V \rightarrow \{-1,+1\}$, and \emph{switching a signed graph $\Sigma = (G,\sigma)$ by $\varsigma$} is the signed graph $\Sigma^{\varsigma}=(G,\sigma^{\varsigma})$ where $\sigma^{\varsigma}=\varsigma(v_i)^{-1}\sigma(e_{ij})\varsigma(v_j)$. The following are well-known facts for signed graphs \cite{SG,OSG}:

\begin{lemma}
Switching is an equivalence relation on the set of signed graphs on an graph. 
\end{lemma}

\begin{lemma}
Switching does not alter the sign of any circle in the signed graph. 
\end{lemma}

\begin{corollary}
The set of balanced signed graphs on a given graph are switching equivalent.
\end{corollary}

\begin{lemma}

Let $\Sigma $ be a signed graph, $[\Sigma ]$ be the switching class of $\Sigma $, and $ne(\Sigma )$ be the number of negative edges in $\Sigma $. The frustration index of signed graph $\Sigma $ is  
\begin{equation*}
fr(\Sigma )=\min\limits_{\Sigma ^{\prime }\in \lbrack \Sigma ]}ne(\Sigma
^{\prime })\text{.}
\end{equation*}
\end{lemma}

For example, all balanced signed graphs have a frustration of $0$ since they can be switched into the all-positive signed graph. 

The classic balancing algorithm is:

\begin{enumerate}
    \item[] \textbf{Signed Graph Balancing Algorithm:}
    \item Input: Signed graph $\Sigma = (G,\sigma)$.
    \item Find a spanning tree $T$ of $G$.
    \item Assign the edges of $T$ the signs in $\Sigma$.
    \item For $e \in \Sigma \setminus T$, assign $e$ the unique sign such that the fundamental cycle is positive.
    \item Output: Balanced signed graph $\Sigma_T$.
\end{enumerate}

Given $\Sigma = (G,\sigma)$, this algorithm determines the nearest balanced signed graphs with respect to a given spanning tree. It seems worthwhile to investigate the set of these nearest balanced signed graphs and their relation to frustration. We extend the idea of balancing sets and frustration to oriented hypergraphs through the underlying incidence structure. Unfortunately, this is not possible for unbalanceable oriented hypergraphs as there is an issue with cross-thetas from Theorem \ref{crossthetabad}. 

\section{Balancing Sets, Shunting, and Balanceable Circuits}

\subsection{Balancing Sets}

Given an oriented hypergraph $G$, a \emph{balancing set of $G$} is the set of incidences whose reversal turns $G$ into a balanced hypergraph. As such, an oriented hypergraph is balanceable if and only if it has a balancing set, thus the underlying hypergraph is cross-theta-free. To find balancing sets of any balanceable oriented hypergraph we translate Camion's Algorithm for re-signing $\{0,1\}$-matrices \cite{Camion} to run on the underlying incidence structure of an oriented hypergraph (see \cite{BM, OHD}).

\begin{enumerate}
    \item[] \textbf{Camion's Signing Algorithm:}
    \item Input: A $\{0,1\}$-matrix $\mathbf{A}$ and its bipartite representation graph $\Gamma$.
    \item Find a spanning tree $T$ of $\Gamma$.
    \item Assign the edges of $T$ arbitrary signs.
    \item For $e \in \Gamma \setminus T$, assign $e$ the unique sign such that the sum of the edge signs of the corresponding fundamental circle is congruent to 0 mod 4.
    \item Output: Balanced matrix $\mathbf{M}$ if $\mathbf{A}$ was balanceable.
\end{enumerate}

We make a trivial adjustment to Camion's Signing Algorithm to apply to the incidences of the oriented hypergraph. Observe that the local circles in oriented hypergraph $G$ are in bijection with the graphic cycles of $\Gamma$, and the 0 mod 4 parity condition from Camion's original algorithm is equivalent to local circles being positive in $G$. Note that this algorithm is also a refinement of the Balancing Algorithm for signed graphs since every edge has exactly $2$ incidences and reorienting a single incidence moves between coherent and introverted/extroverted bidirected edges.

\begin{enumerate}
    \item[] \textbf{Camion's Incidence Reorientation Algorithm:}
    \item Input: An oriented hypergraph $G$ and its oriented bipartite representation graph $\Gamma$.
    \item Find a spanning tree $T$ of $\Gamma$.
    \item Assign the edges of $T$ (incidences of $G$) the signs of the orientation in $G$.
    \item For $e \in \Gamma \setminus T$, assign $e$ the unique sign such that signed of fundamental circle in $\Gamma$ is positive in $G$.
    \item Output: Balanced oriented hypergraph $G'$ if $G$ was balanceable, and a set of incidences $B_T$ whose signs changed.
\end{enumerate}

Directly translating the results of Camion from \cite{Camion} and summarized in \cite{BM} to oriented hypergraphs gives:

\begin{lemma}
An oriented hypergraph is balanceable if, and only if, Camion's Incidence Reorientation Algorithm produces a balanced oriented hypergraph using balancing set $B_T$.
\end{lemma}

\begin{corollary}
Camion's Algorithm produces a balanced oriented hypergraph if, and only if, the original oriented hypergraph was cross-theta-free.
\end{corollary}

\begin{figure}[H]
    \centering
    \includegraphics{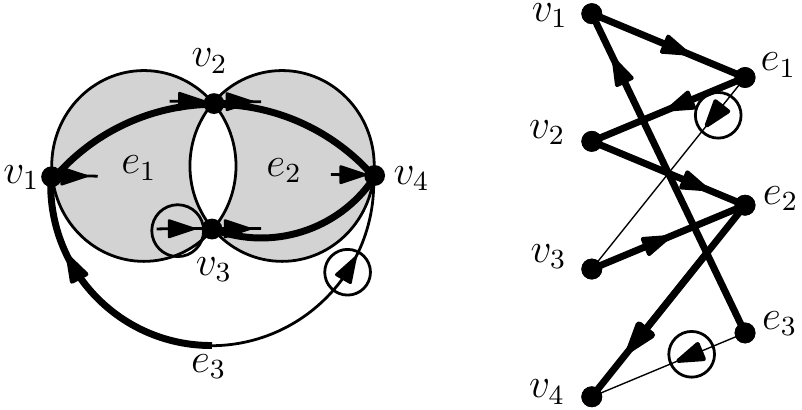}
    \caption{A spanning tree in both $G$ and $\Gamma$ with fundamental incidences/edges circled.}
    \label{fig:Camion}
\end{figure}

Camion proved that the balanced matrix $\mathbf{M}$ produced by his algorithm is unique up to multiplying rows and columns by -1 --- this is equivalent to vertex and edge switching between balancing sets in the oriented hypergraph. 

\begin{lemma}
\label{T:BalSetSwitch}
Given a balancing set $B$, every balancing set of a balanceable oriented hypergraph $G$ is achievable through a finite sequence of vertex and edge switchings in which you add or remove elements to the balancing set. Specifically, if $B$ and $B'$ are balancing sets of an oriented hypergraph $G$ with corresponding binary indicator vectors $\mathbf{b}, \mathbf{b'} \in \mathbb{Z}_2^{I(G)}$, then there exists an incidence bond space vector $\mathbf{s} \in \mathcal{B}(\Gamma_G) \subseteq \mathbf{Z}_2^I$ such that $\mathbf{b} + \mathbf{s} = \mathbf{b'}$.
\end{lemma}

\begin{figure}[H]
    \centering
    \includegraphics{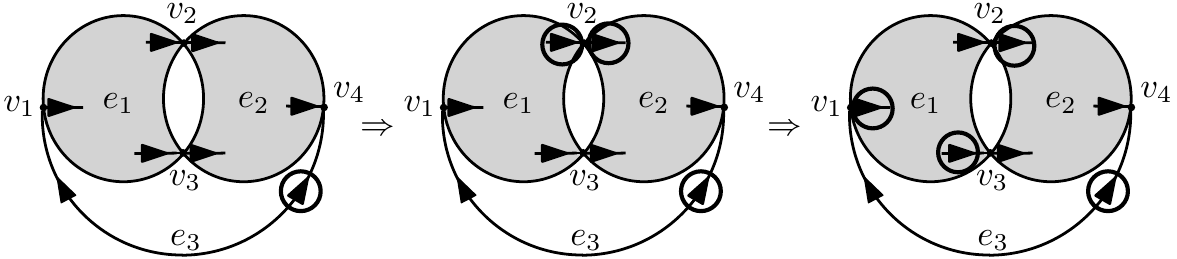}
    \caption{A minimal balancing set (left, circled) and two other balancing sets obtained by switching $v_2$ then $e_1$.}
    \label{fig:camionswitch}
\end{figure}

Additionally, note that the number of edges outside a spanning tree $T$ of $\Gamma$ equals the cyclomatic number of $G$. The resulting balancing set must necessary be minimal.

\begin{lemma}
\label{BNoDisconn}
The deletion of a balancing set disconnects $G$ if, and only if, the balancing set is non-minimal.
\end{lemma}
\begin{proof}
Observe that moving between balanced sets via switching in Lemma \ref{T:BalSetSwitch} is equivalent to adding an element from the binary bond-space of $\Gamma$ to the binary indicator vector of a balanced set. The deletion of a balancing set disconnects $G$ if, and only if, the incidence-bond part of the balancing set of $G$ (edge-bond in $\Gamma$) can be removed by switching by Lemma \ref{T:BalSetSwitch}.
\qed \end{proof}

\begin{lemma}
\label{T:BalSetMin}
Let $G$ be a balanceable hypergraph. $B$ is a minimal balancing set for $G$ if, and only if, $B$ is obtained by Camion’s Incidence Reorientation Algorithm.
\end{lemma}
\begin{proof}
Given a spanning tree $T$ of $\Gamma_G$, let $B_T$ be the balancing set produced by Camion's Incidence Reorientation Algorithm. If $B_T$ is not minimal, there exists a smaller balancing set $S \subset B_T$ and an element $e \in B_T \setminus S$ whose reversal is not necessary to balance $G$. However, by construction, adding $T \cup e$ contains a unique fundamental circle, which must originally be negative prior to the reorientation of $e$, so not reversing $e$ will leave a negative circle. Therefore, $B_T$ must be minimal.

To see the converse, let $B$ be a minimal balancing set for $G$. By Lemma \ref{BNoDisconn} $G \setminus B$ is connected, so any spanning tree of $G \setminus B$ will also be spanning in $G$.
\qed \end{proof}

Given an oriented hypergraph $G$, the minimum number of incidence reversals necessary to balance $G$ is called the \emph{frustration index} of an oriented hypergraph $G$, denoted $fr(G)$. Directly from the previous Lemma we have:

\begin{lemma}
Let $G$ be a balanceable hypergraph and $T$ be a spanning tree of $\Gamma$. Then,
\begin{equation*}
fr(G)=\min\limits_{T}|B_T|
\end{equation*}
where $B_T$ is obtained from Camion’s Incidence Reorientation Algorithm.
\end{lemma}

The main idea is to use incidence reorientation to identify ``weak points'' in oriented hypergraphs to search for additional structure. Unfortunately, this is only well defined on balanceable oriented hypergraphs and an alternate concept seems to be needed for unbalanceable oriented hypergraphs. However, the incidence re-orientation formulation of frustration lines up with acyclic orientations of signed graphs.

\subsection{Shunting}

We introduce an arterial analog of Zaslavsky's \emph{handcuff} characterization of signed graphic circuits containing a negative circle in \cite{SG}. The two critical distinctions are that $1$-edges are treated as balanced pseudo-flowers whose vertex is a thorn, and the single path in a signed graph is replaced with a set of arteries. 

Let $V(B)$ be the (multi-)set of vertices determined by the incidences of a balancing set $B$. Note that $V(B)$ will be a set if $B$ is a minimal balancing set of a balanceable flower or pseudo-flower as every vertex has degree at most $2$. Let $\mathcal{F}$ be a collection of disjoint balanceable flowers and pseudo-flowers where no flower is balanced. Additionally, let $T(\mathcal{F})$ be the set of thorns for each $F \in \mathcal{F}$. A \emph{shunting $\mathcal{S}$ of $\mathcal{F}$} is a collection of disjoint arteries connecting the vertices of balancing set and thorns such that:
\begin{enumerate}
    \item $\mathcal{F} \cup \mathcal{S}$ is connected.
    \item The external vertices of $\mathcal{S}$ are $V(B(\mathcal{F})) \cup T(\mathcal{F})$.
    \item Incidence $i \in B(\mathcal{F})$ if, and only if, there is an $i' \in I(\mathcal{S}$), and their vertices coincide.
\end{enumerate}

An \emph{internal part} of a shunting is any minimal $\mathcal{S}$-path from an $F \in \mathcal{F}$ to itself. An \emph{external part} of a shunting is any minimal $\mathcal{S}$-path between two different elements of $\mathcal{F}$. A shunting is \emph{balanceable} if $\mathcal{F} \cup \mathcal{S}$ is balanceable, hence, cross-theta-free. Internal and external parts of shunts are further refined as follows: a \emph{$tt$-path} is a path between two thorns; a \emph{$bb$-path} is a path between two vertices of a balancing set; and a \emph{$tb$-path (or $bt$-path)} is a path between a thorn and a vertex of a balancing set.

\begin{lemma}
\label{enoughBS}
Let $F$ be a balanceable, but not balanced, flower with distinct vertices $v$ and $w$. If the set of circles that contain $v$ is equal to the set of circles that contain $w$, then $\{v,w\}$ cannot be the vertices of a balancing set.
\end{lemma}

\begin{proof}
Suppose $B = \{i,j\}$ is a balancing set of $F$ with $V(B)=\{v,w\}$.

Take any $vw$-path in the set of $vw$-circles of $F$ and switch $v$ and $w$ as necessary to possibly get a new balancing set $B'=\{i',j'\}$ where $i'$ and $j'$ are in the chosen $vw$-path and $V(B')=\{v,w\}$ still holds. Let $i'$ be incident to edge $e$ and $j'$ be incidence to edge $f$. Since $F$ is a balanceable flower, the degree of each vertex is equal to $2$, and the set of circles containing $v$ are also the set of circles containing $w$, all $vw$-paths within these circles either have first edge $e$ and last edge $f$, or they avoid edges $e$ and $f$ by traversing the circle in the other direction.

Consider the set of $vw$-paths within the $vw$-circles that contain $e$ and $f$. Switch edge $e$ to remove $i'$ from the balancing set and replace it with all of the non-$i'$ incidences of $e$. Next switch the vertices of these new incidences to pass the balancing set to a new set of edges. Continue to switch edges and vertices along these paths as long as all previous switchings have occurred. Since the degree of every vertex is $2$ in a balanceable flower this will terminate with every incidence of $f$ in the balancing set. Switching $f$ provides an empty balancing set, thus $F$ would have to be balanced, a contradiction.
\qed \end{proof}

\begin{lemma}
\label{bscross}
Let $F$ be a balanceable, but not balanced, flower. If $s$ is a single-edge shunt of $F$ corresponding to a balancing set $B$ with $\left\vert B\right\vert \geq 2$, then $F\cup s$ is unbalanceable.
\end{lemma}

\begin{proof}
Let $F$ be a balanceable flower with balancing set $B$ of size at least $2$, and internal shunt $s$.

\textit{Case 1:} If $F$ is a circle-hypergraph consisting of only $2$-edges, then every balancing set contains an odd number of elements. Since $\left\vert B\right\vert \neq 2$, we know that $\left\vert B\right\vert \geq 3$, and the introduction of any internal shunt will produce a cross-theta.

\textit{Case 2:} If $F$ is not a circle-hypergraph, then it must contain an edge of size $3$ or greater.\ Let $v$ and $w$ be different vertices in the vertices of $B$. This can be done since $F$ is a balanceable flower so the degree of every vertex in $F$ is equal to $2$ and there cannot be a balancing set containing only a double incidence since switching out the
double incidence would produce an empty balancing set, making $F$ balanced, a contradiction. 

Let $C$ be a circle in $F$ containing $v$ but not $w$. This can be done by Lemma \ref{enoughBS}. Since $F$ is not a circle-hypergraph there must be an edge $e$ of size $3$ or greater in the circle $C$. Observe that no path from $e$ to $w$ can contain a vertex belonging to circle $C$ or else $F$ would contain a cross-theta, contradicting the fact that $F$ is balanceable. So every $ew$-path is internally disjoint from $C$ or only meets the edges of $C$. If there is an $ew$-path internally disjoint from $C$, then $C$, coupled with internal shunt $s$, forms a cross-theta with end-points $e$ and $v$. 
If there does not exist an $ew$-path internally disjoint from $C$, then take any $ew$-path and let the edge $f\in C$ be the edge closest to $w$. $C$, coupled with internal shunt $s$, form a cross-theta with end-points $f$ and $v$.
\qed \end{proof}

\begin{theorem}
\label{canusemax}
A shunting is balanceable if, and only if, every edge of $\mathcal{S}$ that belongs to a circle is only in $tt$-paths.
\end{theorem}
\begin{proof}
Consider a shunting $\mathcal{F} \cup \mathcal{S}$ and let $s \in E(\mathcal{S})$ belong to a circle in $\mathcal{F} \cup \mathcal{S}$. By construction, $s$ is in no flower-part of any element of $\mathcal{F}$ but belongs to a circle in $\mathcal{F} \cup \mathcal{S}$. Let $F_1,\ldots F_k$ denote the elements of $\mathcal{F}$ that meet any circle containing $s$. Consider the new (pseudo-)flower $F'$ obtained by taking the union of all these $F_i$ along with the elements of each $S \in \mathcal{S}$ connecting them.

If $\mathcal{F} \cup \mathcal{S}$ is balanceable, then by Theorem \ref{ears} the flower-part of $F'$ can be regarded as an ear decomposition consisting of only edge-paths. Again, since $F'$ is balanceable, the degree of each vertex in a circle of $F'$ is equal to $2$, thus $s$ must be in a $tt$-path. Conversely, if $s$ is not in a $tt$-path, then $s$ is either in a $bb$-path or a $tb$-path. In either case the flower-part of $F'$ will have a degree-$3$ vertex, so by Theorem \ref{ohcross}, it must contain a cross-theta.
\qed \end{proof}

\begin{figure}[H]
    \centering
    \includegraphics{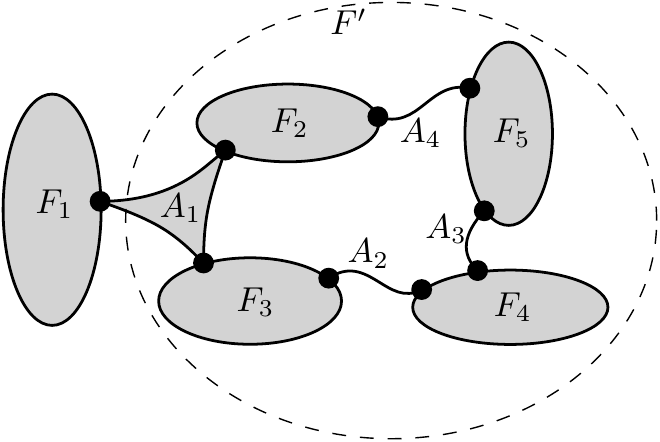}
    \caption{A shunting that forms a larger pseudo-flower.}
    \label{fig:Shunting}
\end{figure}

A shunting $\mathcal{F} \cup \mathcal{S}$ is \emph{$\mathcal{F}$-maximal with respect to $\mathcal{S}$} if, for every non-empty subset $\mathcal{F'} \subseteq \mathcal{F}$ and non-empty edge-induced subhypergraph $\mathcal{S'} \subseteq \mathcal{S}$, $\mathcal{F'} \cup \mathcal{S'}$ is not a flower or pseudo-flower.

\begin{lemma}
\label{usemax}
Let $\mathcal{F} \cup \mathcal{S}$ be $\mathcal{F}$-maximal with respect to $\mathcal{S}$. Every edge $e \in \mathcal{S}$ is an isthmus in $\mathcal{F} \cup \mathcal{S}$.
\end{lemma}
\begin{proof}
From the proof of Theorem \ref{canusemax} if there is a new circle the elements of $\mathcal{F}$ and $\mathcal{S}$ form a new (pseudo-)flower $F'$. The shunting must be balanceable as the adjoining of a cross-path ($tb$-path) on a single balanceable flower makes a larger unbalanced flower.
\qed \end{proof}

\begin{corollary}
Let $\mathcal{F} \cup \mathcal{S}$ be a balanceable $\mathcal{F}$-maximal shunting. The incidence hypergraph $\Upsilon$ with vertex set $\mathcal{F} \cup \mathcal{S}$, edge set $V(B(\mathcal{F})) \cup T(\mathcal{F})$, and incidence set the corresponding shunt incidences, is a tree.
\end{corollary}
\begin{proof}
Immediate from Lemma \ref{usemax} and the fact that every vertex has degree equal to $2$ in a balanceable $\mathcal{F} \cup \mathcal{S}$, so the incidence dual of $V(B(\mathcal{F})) \cup T(\mathcal{F})$ are $2$-edges.
\qed \end{proof}
\begin{figure}[H]
    \centering
    \includegraphics{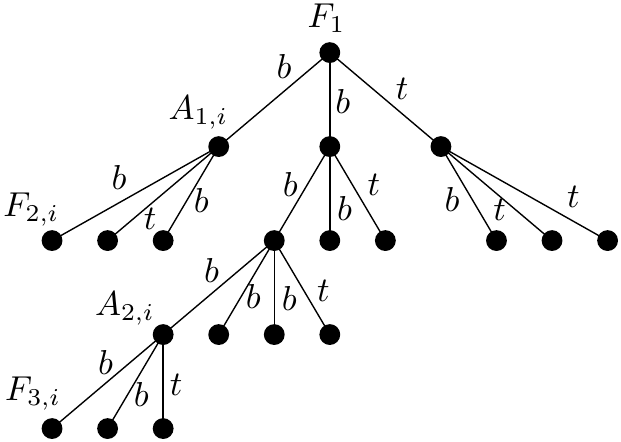}
    \caption{The tree $\Upsilon$ of a balanceable $\mathcal{F}$-maximal shunting.}
    \label{fig:PAbtTree}
\end{figure}

A shunting $\mathcal{F} \cup \mathcal{S}$ is \emph{$\mathcal{S}$-minimal with respect to $\mathcal{F}$} if, for every non-empty subset $\mathcal{F'} \subseteq \mathcal{F}$ and non-empty edge-induced subhypergraph $\mathcal{S'} \subseteq \mathcal{S}$, $\mathcal{F'} \cup \mathcal{S'}$ is not a shunting (for any balancing set).

\begin{lemma}
A shunting $\mathcal{F} \cup \mathcal{S}$ is $\mathcal{S}$-minimal if, and only if, it arises from a minimal balancing set.
\end{lemma}
\begin{proof}
Clearly a minimal balancing set is $\mathcal{S}$-minimal.

To see the other direction, let $N$ be a non-minimal balancing set that forms a shunting. From Lemma \ref{BNoDisconn} there must be a bond in some $F \in \mathcal{F}$ that can be removed from $N$ via switching (Lemma \ref{T:BalSetSwitch}). So the shunting resulting from $N$ is not $\mathcal{S}$-minimal.
\qed \end{proof}

A shunting $\mathcal{F} \cup \mathcal{S}$ that is both $\mathcal{F}$-maximal with respect to $\mathcal{S}$ and $\mathcal{S}$-minimal with respect to $\mathcal{F}$ is called an \emph{optimal shunting}.
 
\begin{lemma}
\label{SGShuntIsDep}
$\mathcal{F} \cup \mathcal{S}$ is an optimal shunting of a signed graph if, and only if, $\mathcal{F} \cup \mathcal{S}$ is a circuit of the signed graphic frame matroid that is not a positive-circle.
\end{lemma}
\begin{proof}
From Lemmas \ref{onlyflowers} and \ref{onlyartery} the only flowers are circle subgraphs and the only arteries are paths, while the only pseudo-flowers of a signed graph are $1$-edges. Since $\mathcal{F} \cup \mathcal{S}$ is a signed graph, from Lemma \ref{bscross} all minimal balancing sets have size equal to $1$. Thus, an optimal shunting in a signed graph consists of two negative circles connected by a path, where the negative circles may be replaced with $1$-edges --- which are the non-positive-circle circuits of the signed graphic frame matroid.
\qed \end{proof}

\subsection{Balanceable Circuits}

With Lemma \ref{SGShuntIsDep} we see that optimal shunting completes the characterization of signed graphic circuits. We now prove that the shunting construction produces balanceable oriented hypergraphic circuits. However, as in \cite{OH1}, the arteries of a balanceable shunting can be $2$-vertex contracted so that the thorns of the pseudo-flowers are removed and connected pseudo-flowers share a common edge. A \emph{$t,k$-hypercircle} is a hypergraph with $t$ monovalent vertices such that balanced subdivision produces an arterial connection of $k$ pseudo-flowers by thorns. Alternatively, a $t,k$-hypercircle is the $2$-vertex contraction of the vertices of an arterial connection. A loose edge is a $0,0$-hypercircle and a non-loose-edge flower is $0,1$-hypercircle.

The introduction of hypercircles in \cite{OH1} is done to provide a direct characterization of the column dependencies of an incidence matrix $\mathbf{H}_G$ that correspond to an edge-induced subhypergraph of $G$ and prevent over-use of Lemma \ref{balsubgood}. By construction, shunting extends to hypercircles. Let $\mathcal{H}$ be a collection of nearly disjoint hypercircles. A \emph{shunting $\mathcal{S}$ of $\mathcal{H}$} is a shunting on flower and pseudo-flower-parts of $\mathcal{H}$. The remaining shunting definitions are defined similarly.

\begin{theorem}
Let $G$ be a balanceable oriented hypergraph with incidence matrix $\mathbf{H}_G$. $\mathbf{H}_G$ is minimally dependent if, and only if, $G$ is a balanced subdivision of an optimal shunting of hypercircles.
\end{theorem}

\begin{proof}
By Lemma \ref{balsubgood} we only need to show optimal shunting of hypercircles is minimally dependent.

Let $\mathcal{H} \cup \mathcal{S}$ be an optimal shunting of hypercircles corresponding to balancing set $\mathcal{B}$, where the elements of $\mathcal{H}$ are $m$ balanceable $t_i,k_i$-hypercircles $H_i$ for $ 1 \leq i \leq m $, and the flower-parts of each $H_i$ are $F_{i,j}$ for $1 \leq j \leq k_i$. 

By Lemma \ref{usemax} if $C$ is a circle in $\mathcal{H} \cup \mathcal{S}$, then $C$ belongs to some $H_i$, hence, some $F_{i,j}$ by Lemma \ref{balsubgood}. So the cyclomatic number of $\mathcal{H} \cup \mathcal{S}$ is determined by the individual flower-parts:
\begin{align}
\label{E:Cyclo}
    \varphi_{\mathcal{H} \cup \mathcal{S}} = \sum\limits_{i=1}^{m}\varphi
_{H_{i}} = \sum\limits_{i=1}^{m} \sum\limits_{j=1}^{k_i} \varphi
_{F_{i,j}}.
\end{align}

By Lemma \ref{T:BalSetMin} there is a spanning tree $T$ of $\Gamma_{\mathcal{H} \cup \mathcal{S}}$ that produces $\mathcal{B}$. Since the only circles belong to the flower-parts of the hypercircles we can partition $\mathcal{B}$ into balancing sets $B_{H_1} , B_{H_2} , \ldots , B_{H_m}$, which can be further partitioned into balancing sets $B_{i,j}$ for each $F_{i,j}$.
\begin{align}
\label{E:NegCirc}
    \vert \mathcal{B} \vert = \sum\limits_{i=1}^{m} \vert B_{H_{i}} \vert = \sum\limits_{i=1}^{m} \sum\limits_{j=1}^{k_i} \vert B_{F_{i,j}} \vert.
\end{align}

For each $F_{i,j}$, take a system of distinct representatives of the vertices of the positive essential circles along the associated edge-ear decomposition from Lemma \ref{ears}. For each vertex in the system of distinct representatives, there is a linear combination of rows of $\mathbf{H}_{\mathcal{H} \cup \mathcal{S}}$ that produces a row of zeroes in the corresponding square sub-matrix of the positive locally-signed-graphic essential circle. However, by Lemmas \ref{ohcross} and \ref{crossthetabad} each of these vertices has degree equal to $2$ in $\mathcal{H} \cup \mathcal{S}$ since they are in a balanceable flower-part and are not a vertex of their respective balancing set. Thus, the row combinations produce an entire row of zeroes in $\mathbf{H}_{\mathcal{H} \cup \mathcal{S}}$. The columns of sub-matrices for negative locally-signed-graphic circles are independent, so the row rank of $\mathbf{H}_{\mathcal{H} \cup \mathcal{S}}$ is
\begin{align*}
r_{ \mathcal{H} \cup \mathcal{S}} &= \vert V_{\mathcal{H} \cup \mathcal{S}} \vert - p.
\end{align*}
where $p$ is the number of positive essential circles in $\mathcal{H} \cup \mathcal{S}$.

From Equations \ref{E:Cyclo} and \ref{E:NegCirc} there are 
\begin{align*}
\varphi_{\mathcal{H} \cup \mathcal{S}} - \sum\limits_{i=1}^{m} \sum\limits_{j=1}^{k_i} \vert B_{F_{i,j}} \vert = \varphi_{\mathcal{H} \cup \mathcal{S}} - \sum\limits_{i=1}^{m} \vert B_{H_{i}} \vert = \varphi_{\mathcal{H} \cup \mathcal{S}} - \vert \mathcal{B} \vert 
\end{align*}
positive essential circles so the row rank is
\begin{align*}
r_{ \mathcal{H} \cup \mathcal{S}} &= \vert V_{\mathcal{H} \cup \mathcal{S}} \vert - \left(\varphi_{\mathcal{H} \cup \mathcal{S}} - \vert \mathcal{B} \vert \right).
\end{align*}

In order for $\mathbf{H}_{\mathcal{H} \cup \mathcal{S}}$ to be minimally dependent, the nullity of $\mathbf{H}_{\mathcal{H} \cup \mathcal{S}}$ must be equal to $1$, and no edge-induced subhypergraph can be dependent. Using Equations \ref{E:Cyclo} and \ref{E:NegCirc} the cyclomatic number of $\mathcal{H} \cup \mathcal{S}$ is
\begin{align*}
\varphi_{\mathcal{H} \cup \mathcal{S}} &= \vert I_{\mathcal{H} \cup \mathcal{S}} \vert - (\vert V_{\mathcal{H} \cup \mathcal{S}} \vert + \vert E_{\mathcal{H} \cup \mathcal{S}} \vert) + 1 \\
&=  2 \vert V_{\mathcal{H} \cup \mathcal{S}} \vert + \vert \mathcal{B} \vert - (\vert V_{\mathcal{H} \cup \mathcal{S}} \vert + \vert E_{\mathcal{H} \cup \mathcal{S}} \vert) + 1 \\
&= \vert V_{\mathcal{H} \cup \mathcal{S}} \vert + \vert \mathcal{B} \vert -  \vert E_{\mathcal{H} \cup \mathcal{S}} \vert + 1.
\end{align*}

Solving for $\vert E_{\mathcal{H} \cup \mathcal{S}} \vert - 1$ we have
\begin{align*}
\vert E_{\mathcal{H} \cup \mathcal{S}} \vert - 1 &= \vert V_{\mathcal{H} \cup \mathcal{S}} \vert - \left( \varphi_{\mathcal{H} \cup \mathcal{S}} - \vert \mathcal{B} \vert \right) \\
&= r_{ \mathcal{H} \cup \mathcal{S}}.
\end{align*}
Thus, $\mathbf{H}_{\mathcal{H} \cup \mathcal{S}}$ is a nullity-$1$ matrix. Moreover, every proper edge-induced subgraph either contains a monovalent vertex or leaves an unshunted vertex of a balancing set. Since we started with an optimal shunting, $\mathcal{F} \cup \mathcal{S}$ is minimally dependent.

Now assume that $G$ is not an optimal shunting of hypercircles. If $G$ is disconnected, or contains a monovalent vertex, it cannot be minimally dependent. Therefore, $G$ must be a cross-theta-free oriented hypergraph where every vertex has degree equal to $2$ or greater. We may also assume that $G$ has a negative circle $C$, or else Theorem \ref{BalMD} applies. Thus, $C$ must be contained in some flower-part of some hypercircle $H$ whose with a non-empty balancing set. If $G = H$ it is not minimally dependent. If $G \neq H$, then $G \setminus H$ is non-empty, and if $G$ contains any part of a non-optimal shunting then it either properly contains an optimal shunting (hence, contain a minimal dependency), or by Lemma \ref{bscross} would be unbalanceable. Finally, if a part of $G$ avoids a shunting entirely there there is some negative circle that does not connect to a shunt, and it cannot be minimally dependent.\qed 
\end{proof}

\begin{example}
A simple example tells us how to adjoin elementary basis vectors as columns to the incidence matrix so that a circuit is formed. Let $\mathcal{F} \cup \mathcal{S}$ be a shunting with minimal balancing set $B$ with $V(B) = \{v_1, \ldots v_{\vert B \vert}\}$. Let $\mathcal{F} = \{F , P_1 , \ldots , P_{\vert B \vert}$\} where $F$ is a balanceable flower, each $P_i$ is a $1$-edge pseudo-flower $\{v_i , e_i \}$, and $\mathcal{S} = \{v_1, \ldots v_{\vert B \vert}\}$. The shunting $\mathcal{F} \cup \mathcal{S}$ is optimal, hence, minimally dependent.
\end{example}

\section{A Note on Unbalanceable Circuits}

\subsection{Minimal cross-thetas}

The cross-theta plays a central role in completing the circuit characterization of the oriented hypergraphic matroid. The techniques so far rely on either being balanced or balanceable with Camion's algorithm to find balancing sets of incidences and providing a reinterpretation of frustration, all of which require cross-theta-free hypergraphs. An alternative formulation of ``frustration'' seems necessary to tackle unbalanceable oriented hypergraphs. We examine some simple properties of minimal cross-thetas to provide further context in the importance and difficulty of unbalanced oriented hypergraphs.

As discussed in \cite{OH1} the proof techniques apply to any matrix whose entries are commensurable. The unit element is represented by an entrant arrow and other entries are represented by multiple arrows. An entry of $\pm 3$ in an incidence matrix, represented as $3$ entrant (or salient) arrows, is the smallest cross-theta. In fact, reorienting one of these arrows produces a value of $\pm 1$, which provides the missing hypergraph family from the characterization of totally unimodular matrices as discussed in \cite{BM}.

Let $L_k$ denote the hypergraph consisting of a single vertex, a single edge, and $k$ incidences. $L_k$ is \emph{extroverted} if each incidence is $+1$, and \emph{introverted} if each incidence is $-1$.  A $k$-cross-theta ($k \geq 3$) is a subhypergraph that consists of $k$ internally disjoint paths of half-integer length.

\begin{lemma}
Every minimal $k$-cross-theta is a subdivision of $L_k$.
\end{lemma}
\begin{proof}
$L_k$ consists of $k$ paths of length $1/2$. Subdivision increases the length of a path by integer length.
\qed \end{proof}

The following corollaries are immediate from the observation that balanced subdivision preserves minimal dependency.

\begin{corollary}
Every circle in a minimal $k$-cross-theta is negative if, and only if, it is switching equivalent to a balanced subdivision of an extroverted or introverted $L_k$.
\end{corollary}

\begin{corollary}
A minimal $k$-cross-theta in which every circle is negative is minimally dependent over $GF(k)$.
\end{corollary}

\begin{corollary}
A minimal $(p+n)$-cross-theta that is a balanced subdivision of an $L_{p+n}$ with $p$ entrant and $n$ salient arrows is minimally dependent over $GF(\left\vert p-n \right\vert)$.
\end{corollary}

\begin{corollary}
A minimal $2k$-cross-theta that is a balanced subdivision of an $L_{2k}$ with $k$ entrant and $k$ salient arrows is minimally dependent over every field.
\end{corollary}

It was discussed in \cite{OH1} that $3$-cross-thetas must have a negative circle. Since every $k$-cross-theta contains a $3$-cross-theta they all trivially must contain a negative circle. However, the precise minimum number of negative circles in a minimal $k$-cross-theta is given by the following Lemma. The minimum number of negative circles possible over all orientations may provide the appropriate alternative for frustration.

\begin{lemma}
\label{quartsqare}
The minimum number of negative circles in a minimal $k$-cross-theta is
$\dbinom{\left\lfloor \frac{k}{2}\right\rfloor }{2}+\dbinom{\left\lfloor \frac{k+1}{2}\right\rfloor }{2} =\left\lfloor \frac{k-1}{2}\right\rfloor \left\lfloor \frac{k}{2} \right\rfloor =  \left\lfloor \frac{k-1}{2} \right\rfloor \left\lceil \frac{k-1}{2}\right\rceil = \left\lfloor \frac{(k-1)^{2}}{4}\right\rfloor$.
\end{lemma}
\begin{proof}
The maximum number of negative circles occur at an extroverted or introverted $L_k$. The minimum occurs when half are entrant/salient. The remaining equalities are similar expressions of the quarter-squares sequence.
\qed \end{proof}

\begin{theorem}
If $F$ is a minimal $k$-cross-theta that does not vanish over $GF(q)$, and $P$ is a $1$-edge pseudo-flower that shares its vertex with $F$, then $F \cup P$ is minimally dependent.
\end{theorem}
\begin{proof}
$F \cup P$ has $\left\vert E_{F \cup P} \right\vert = \left\vert V_{F \cup P} \right\vert + 1$, and $F$ does not vanish over $GF(q)$.
\qed \end{proof}

\subsection{Example: The Fano and non-Fano Matroids}

\begin{figure}[H]
    \centering
    \includegraphics[scale=1.0]{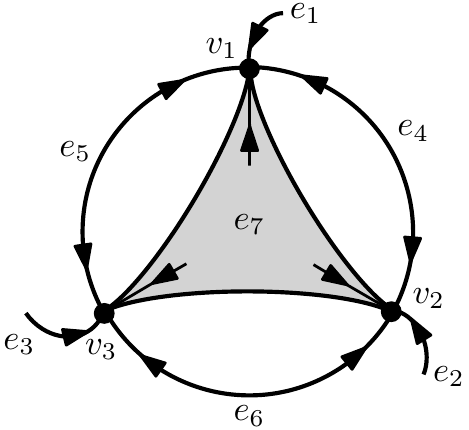}
    \caption{An extroverted complete hypergraph on $3$ vertices.}
    \label{fig:Fano}
\end{figure}

A \emph{complete hypergraph} is a hypergraph whose edges correspond to the faces of a simplex. The incidence matrix 
\begin{align*}
\mathbf{H} =
\left[
\begin{array}{ccccccc}
1 & 0 & 0 & 1 & 1 & 0 & 1 \\ 
0 & 1 & 0 & 1 & 0 & 1 & 1 \\ 
0 & 0 & 1 & 0 & 1 & 1 & 1 
\end{array}
\right]
\end{align*}
has its corresponding complete hypergraph depicted in Figure \ref{fig:Fano}.

\begin{figure}[H]
    \centering
    \includegraphics[scale=1.0]{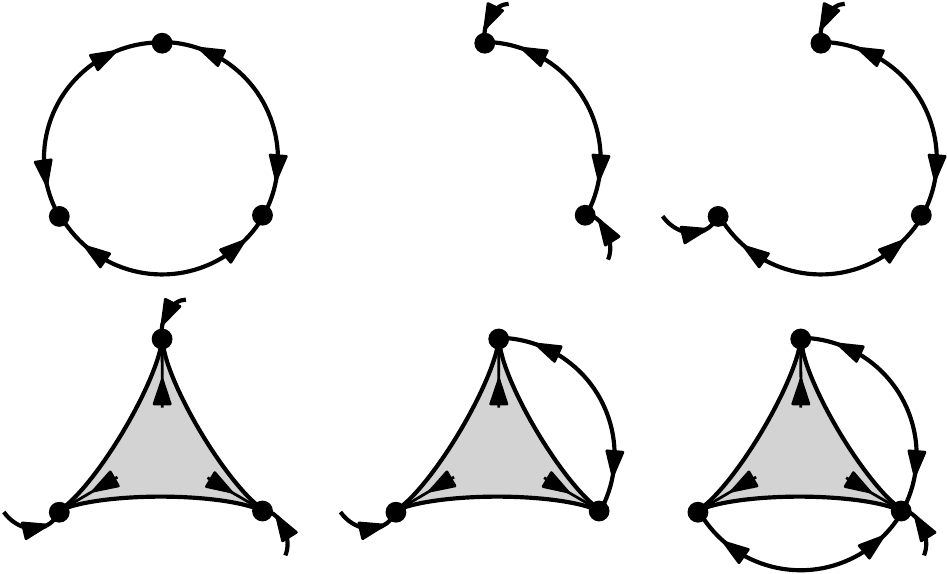}
    \caption{The circuit families of the Fano matroid.}
    \label{fig:FanoCirc}
\end{figure}

Since every edge is extroverted in Figure \ref{fig:Fano}, every adjacency is negative, and the only unbalanced circle is $C = \{v_1, e_4, v_2, e_6, v_3, e_5, v_1\}$. However, if we regard $\mathbf{H}$ as a matrix over $GF(2)$ there are no negative circles so circle $C$ is balanced and minimally dependent, while the remainder of the shunting families are always minimally dependent. Thus, the hypergraph families in Figure \ref{fig:FanoCirc} are all the circuits of the Fano matroid.

However, if we regard $\mathbf{H}$ as a matrix over $GF(3)$ circle $C$ immediately switches from balanced to unbalanced --- which is easily recognized as relaxing the circuit hyperplane to produce the non-Fano matroid. We must then shunt $C$ by any of the $1$-edge shunts or the $3$-edge shunt via a non-minimal-balancing-set as described in Lemma \ref{bscross}.

\begin{figure}[H]
    \centering
    \includegraphics[scale=1.0]{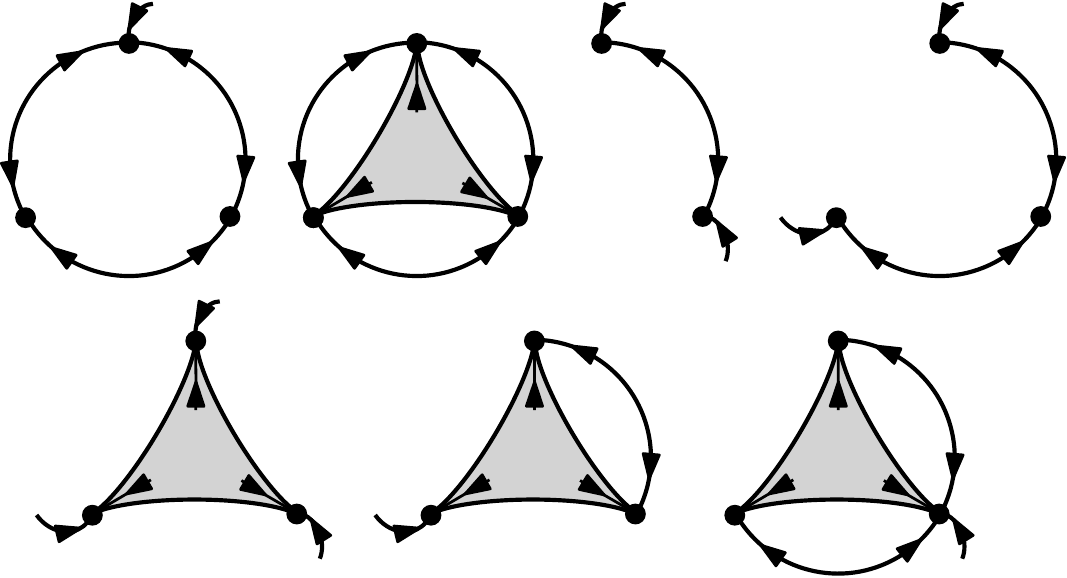}
    \caption{The circuit families of the Non-Fano matroid.}
    \label{fig:NonFanoCirc}
\end{figure}

With the incorporation of the unbalancable circuits the hope is to then apply locally signed-graphic techniques to representable matroids to provide greater understanding of the connection between graphs and matroids. Moreover, with the direct connection from oriented hypergraphs to the bipartite incidence graph $\Gamma$ (via a logical functor from \cite{IH1}) it opens the door for a modified version of graph theoretic techniques, such as Robertson-Seymour, to be applied to representable matroids.

%%%%%%%%%%%%%%%%%%%%%%%%%%%%%%%%
\newpage
\section*{References}
\bibliographystyle{amsplain2}
\bibliography{mybibA}

\end{document}